\documentclass[12pt]{article}

\usepackage{amssymb,amsmath,amsfonts,enumerate}

\title{
\author{\bf{Lech Pasicki}}
\bf{My last fixed point theorem}
}

\newtheorem{theorem}{\indent Theorem}

\newtheorem{definition}[theorem]{\indent Definition}

\newcommand{\n}{\mbox{$_{n}$}}

\newcommand{\no}{\mbox{$n_{0}$}}
\newcommand{\nl}{\mbox{$_{n+1}$}}
\newcommand{\nn}{\mbox{$n \in \mathbb{N}$}}
\newcommand{\xnn}{\mbox{$(x_{n})_{n \in \mathbb{N}}$}}
\newcommand{\pxnm}{\mbox{$p(x_{n},x_{m})$}}
\newcommand{\pxnx}{\mbox{$p(x_{n},x)$}}
\newcommand{\limn}{\mbox{$\lim_{m,n \to \infty}$}}

\newcommand{\li}{\mbox{$\lim_{n \to \infty}$}}

\newcommand{\xyX}{\mbox{$x,y \in X$}}

\newcommand{\fno}{\mbox{$f^{n}x_{0}$}}

\begin{document}
\maketitle
\vspace{1 in}

\begin{abstract}
In this note a far extension of the Banach fixed point theorem is proved.
\end{abstract}
%a.1
\begin{definition}[{\rm \cite[Definition 6.1]{p40}\bf}]
\label{De1}
A mapping $p\colon X \times X \to [0,\infty)$ is a felt metric for $X$ if the following system of conditions 
is satisfied: 
\begin{subequations}
\label{con1}
\begin{align}
 &p(x,y) = 0 \textit{ yields } x = y, \quad \xyX, \label{con1a}\\
 &p(x,y) = p(y,x), \quad \xyX, \label{con1b}\\
\begin{split}
&\textit{for each } \epsilon > 0 \textit{ there exists a } \delta > 0 \textit{ such that }\\
&p(z,y) < \delta \textit{ yields } |p(z,x) - p(y,x)| < \epsilon, \quad x,y,z, \in X. \label{con1c}
\end{split}
\end{align}
\end{subequations}
A felt metric space $(X,p)$ (its topology is generated by balls) is $0$-complete, if for each sequence 
$\xnn$ in $X$, such that $\limn \pxnm = 0$, there exists an $x \in X$ for which $\li \pxnx = 0$; a selfmapping $f$ 
on $X$ is $0$-continuous at an $x \in X$, if $\li \pxnx = 0$ yields $\li p(fx\n,fx) = 0$, for each sequence 
$\xnn$ in $X$; $f$ is $0$-continuous if it is $0$-continuous at each $x \in X$.
\end{definition}
\par In \cite[Theorem 6.6]{p40} the following condition was used:
\begin{equation}
\label{con2}
\begin{split}
  &\textit{for each } \alpha > 0 \textit{ there exists an } \epsilon > 0 \textit{ such that}\\
  &\alpha - \epsilon < p(y,x) < \alpha + \epsilon  \textit{ yields } p(fy,fx) \leq \alpha, \quad \xyX. 
\end{split}
\end{equation}
\par Now, let us prove a result stronger than \cite[Theorem 6.6]{p40}.
%a.4
\begin{theorem}
\label{a.4}
Assume that $p$ is a felt metric for $X$, $f$ is a selfmapping on $X$, $x\n = \fno$, $\nn$, $\li p(x\nl,x\n) = 0$, 
and let \eqref{con2} or 
\begin{equation}
\label{con3}
\begin{split}
  &\textit{for each } \alpha > 0 \textit{ there exists an } \epsilon > 0 \textit{ such that}\\
  &\alpha \leq p(y,x) < \alpha + \epsilon  \textit{ yields } p(fy,fx) \leq \alpha, \quad \xyX 
\end{split}
\end{equation}
be satisfied (equivalent conditions). If $(X,p)$ is $0$-complete, then there exists an $x \in X$ such that 
$\li p(x\n,x) = p(x,x) = 0$ and $fx = x$. 
\end{theorem}
Proof. In particular, for $\alpha = p(y,x) > 0$, condition \eqref{con2} as well as \eqref{con3} yields 
$p(fy,fx) \leq p(y,x)$, {\it i.e.} $f$ is nonexpansive (for $p(y,x) > 0$), and conditions \eqref{con2}, 
\eqref{con3} are equivalent (one may assume $\alpha - \epsilon > 0$ in \eqref{con2}). Now, in view of 
\cite[Theorem 6.6]{p40}, if $(X,p)$ is $0$-complete, then for some $x \in X$ we have 
$\li p(x\n,x) = p(x,x) = 0$. From \eqref{con1c} it follows that
\begin{equation}
\label{con4}
 \li |p(x\nl,fx) - p(x,fx)| = \li p(x\nl,x) = 0.
\end{equation}
If $p(x\n,x) = 0$ holds for all $n > \no$, then $x = x\n = x\nl$ and we get $fx = fx\n = x\nl = x$. Therefore, 
one may assume $p(x\n,x) > 0$ for infinitely many $\nn$. Suppose $p(x,fx) = 2\beta > 0$. Then 
(see \eqref{con4}) we obtain $p(x\nl,fx) > \beta$ for large $\nn$. Let us consider an $\nn$ such that 
$0 < p(x\n,x) = \alpha < \beta$. Now, \eqref{con3} yields  $\beta < p(x\nl,fx) \leq \alpha$, a 
contradiction. Consequently, $p(x,fx) = 0$, and $fx = x$ (see \eqref{con1a}).$\quad \square$
\par Theorem \ref{a.4} is a far extension of the Banach theorem, as for $0 \leq c < 1$ one can calculate the 
respective $\epsilon > 0$ (see \eqref{con3}) from the subsequent formula 
\[
 p(fy,fx) \leq cp(y,x) \leq c(\alpha + \epsilon) \leq \alpha.
\]
\section*{Acknowledgements}
This work was partially supported by the Faculty of Applied Mathematics AGH UST statutory tasks within subsidy 
of the Polish Ministry of Science and Higher Education, grant no. 16.16.420.054.
\par 
%\vspace{.1in}
\mbox{Faculty of Applied Mathematics} \linebreak
\mbox{AGH University of Science and Technology} \linebreak
\mbox{Al. Mickiewicza 30} \linebreak
\mbox{30-059 KRAK\'OW, POLAND} \linebreak
\mbox{E-mail: pasicki@agh.edu.pl}


\vspace{.5in}
\begin{thebibliography}{0}
\bibitem{p40}
L. Pasicki, A strong fixed point theorem, Topology Appl., 282 (2020) 107300, DOI: 10.1016/j.topol.2020.107300.
\end{thebibliography}
\end{document}